\documentclass[12pt]{amsart}
\input{epsf}
\usepackage{amsmath,amsfonts,amssymb,latexsym,fullpage}

\usepackage{hyperref}

\newtheorem{theorem}{Theorem}[section]

\newtheorem{lemma}[theorem]{Lemma}
\theoremstyle{definition}

\newtheorem{corollary}[theorem]{Corollary}

\theoremstyle{remark}
\newtheorem{remark}[theorem]{Remark}
\numberwithin{equation}{section}

\def\D{\mathfrak D}
\def\DD{\mathcal D}

\def\f{f_1}
\def\ff{f_2}
\def\kra{\varphi}
\def\kth{\psi}
\def\inj{g_1}
\def\injj{g_2}
\def\ul{\underline}
\def\lis{\mathrm{lis}}
\def\lds{\mathrm{lds}}
\renewcommand\S{\mathfrak{S}}
\renewcommand\u{{\mathbf u}}
\renewcommand\d{{\mathbf d}}

\author[A. Burstein]{Alexander Burstein}
\address{Department of Mathematics,
Iowa State University, Ames, IA 50011-2064 USA}
\email{burstein@math.iastate.edu}

\author[S. Elizalde]{Sergi Elizalde}
\address{Department of Mathematics, Dartmouth College,
Hanover, NH 03755-3551 USA} \email{sergi.elizalde@dartmouth.edu}

\author[T. Mansour]{Toufik Mansour}
\address{Department of Mathematics, University of Haifa,
31905 Haifa, Israel}
\email{toufik@math.haifa.ac.il}

\title[Restricted Dumont permutations]{Restricted Dumont permutations, Dyck
paths, and noncrossing partitions}



\begin{document}

\begin{abstract}
We complete the enumeration of Dumont permutations of the second
kind avoiding a pattern of length 4 which is itself a Dumont
permutation of the second kind. We also consider some combinatorial
statistics on Dumont permutations avoiding certain patterns of
length 3 and 4 and give a natural bijection between $3142$-avoiding
Dumont permutations of the second kind and noncrossing partitions
that uses cycle decomposition, as well as bijections between 132-,
231- and 321-avoiding Dumont permutations and Dyck paths. Finally,
we enumerate Dumont permutations of the first kind simultaneously
avoiding certain pairs of 4-letter patterns and another pattern of
arbitrary length.
\end{abstract}

\maketitle  

\section{Preliminaries}\label{sec:prelim}
The main goal of this paper is to give analogues of enumerative
results on certain classes of permutations characterized by
pattern-avoidance in the symmetric group $\S_n$. In the set of
Dumont permutations (see below) we identify classes of restricted
Dumont permutations with enumerative properties analogous to results
on permutations. More precisely, we study the number of Dumont
permutations of length $2n$ avoiding either a 3-letter pattern or a
4-letter pattern. We also give direct bijections between
equinumerous sets of restricted Dumont permutations of length $2n$
and other objects such as restricted permutations of length $n$,
Dyck paths of semilength $n$, or noncrossing partitions of
$[n]=\{1,2\dots,n\}$.

\subsection{Patterns}\label{subsec:patterns}
Let $\sigma\in\S_n$ and $\tau\in\S_k$ be two permutations. We say
that $\tau$ \emph{occurs} in $\sigma$, or $\sigma\in\S_n$
\emph{contains} $\tau$, if $\sigma$ has a subsequence
$(\sigma(i_1),\dots,\sigma(i_k))$, $1\le i_1<\dots<i_k\le n$, that
is order-isomorphic to $\tau$ (in other words, for any $j_1$ and
$j_2$, $\sigma(i_{j_1})\le\sigma(i_{j_2})$ if and only if
$\tau(j_1)\le\tau(j_2)$). Such a subsequence is called an
\emph{occurrence} (or an \emph{instance}) of $\tau$ in $\sigma$. In
this context, the permutation $\tau$ is called a \emph{pattern}. If
$\tau$ does not occur in $\sigma$, we say that $\sigma$
\emph{avoids} $\tau$, or is \emph{$\tau$-avoiding}. We denote by
$\S_n(\tau)$ the set of permutations in $\S_n$ avoiding a pattern
$\tau$. If $T$ is a set of patterns, then $\S_n(T)=\cap_{\tau\in
T}{\S_n(\tau)}$, i.e. $\S_n(T)$ is the set of permutations in $\S_n$
avoiding all patterns in $T$.

The first results in the extensive body of research on permutations
avoiding a 3-letter pattern are due to Knuth \cite{Knuth}, but the
intensive study of patterns in permutations began with Simion and
Schmidt \cite{SS} who considered permutations and involutions
avoiding each set $T$ of 3-letter patterns. One of the most
frequently considered problems is the enumeration of $\S_n(\tau)$
and $\S_n(T)$ for various patterns $\tau$ and sets of patterns $T$.
The inventory of cardinalities of $|\S_n(T)|$ for $T\subseteq\S_3$
is given in \cite{SS}, and a similar inventory for
$|\S_n(\tau_1,\tau_2)|$, where $\tau_1\in\S_3$ and $\tau_2\in\S_4$
is given in \cite{West2}. Some results on $|\S_n(\tau_1,\tau_2)|$
for $\tau_1,\tau_2\in\S_4$ are obtained in \cite{West1}. The exact
formula for $|\S_n(1234)|$ and the generating function for
$|\S_n(12\dots k)|$ are found in \cite{Gessel}. B\'{o}na \cite{Bona}
has found the exact value of $|\S_n(1342)|=|\S_n(1423)|$, and
Stankova \cite{Stankova1,Stankova2} showed that
$|\S_n(3142)|=|\S_n(1342)|$. For a survey of results on pattern
avoidance, see \cite{Bona-book, KM}.

Another problem is finding equinumerously avoided (sets of)
patterns, i.e. sets $T_1$ and $T_2$ such that
$|\S_n(T_1)|=|\S_n(T_2)|$ for any $n\ge 0$. Such (sets of) patterns
are called \emph{Wilf-equivalent} and said to belong to the same
\emph{Wilf class}. There are eight symmetry operations on $\S_n$
that map every pattern onto a Wilf-equivalent
pattern, including:
\begin{itemize}

\item \emph{reversal} $r$: $r(\tau)(j)=\tau(n+1-j)$, i.e.
$r(\tau)$ is $\tau$ read right-to-left.

\item \emph{complement} $c$: $c(\tau)(j)=n+1-\tau(j)$, i.e.
$c(\tau)$ is $\tau$ read upside down.

\item $r\circ c=c\circ r$: $r\circ c(\tau)(j)=n+1-\tau(n+1-j)$,
i.e. $r\circ c(\tau)$ is $\tau$ read right-to-left upside down.

\item \emph{inverse} $i$: $i(\tau)=\tau^{-1}$.
\end{itemize}
The set of patterns $\langle r,c,i\rangle (\tau)
=\{\tau,r(\tau),c(\tau),r(c(\tau)),\tau^{-1},r(\tau^{-1}),c(\tau^{-1}),r(c(\tau^{-1}))\}$
is called the \emph{symmetry class} of $\tau$.

Sometimes we will represent a permutation $\pi\in\S_n$ by placing
dots on an $n\times n$ board. For each $i=1,\ldots,n$, we will place
a dot with abscissa $i$ and ordinate $\pi(i)$ (the origin of the
board is at the bottom-left corner).

\subsection{Dumont permutations}\label{subsec:dumont}
In this paper we answer some of the above problems in the case of
Dumont permutations. A \emph{Dumont permutation of the first kind}
is a permutation $\pi\in\S_{2n}$ where each even entry is followed
by a descent and each odd entry is followed by an ascent or ends the
string. In other words, for every $i=1,2,\dots,2n$,
\[
\begin{split}
\pi(i) \text{ is even} &\implies i<2n \text{ and } \pi(i)>\pi(i+1),\\
\pi(i) \text{ is odd}  &\implies \pi(i)<\pi(i+1) \text{ or } i=2n.
\end{split}
\]
A \emph{Dumont permutation of the second kind} is a permutation
$\pi\in\S_{2n}$ where all entries at even positions are deficiencies
and all entries at odd positions are fixed points or excedances. In
other words, for every $i=1,2,\dots,n$,
\[
\begin{split}
\pi(2i)&<2i,\\
\pi(2i-1)&\ge 2i-1.
\end{split}
\]
We denote the set of Dumont permutations of the first (resp. second)
kind of length $2n$ by $\D^1_{2n}$ (resp. $\D^2_{2n}$). For example,
$\D^1_{2}=\D^2_{2}=\{21\}$, $\D^1_{4}=\{2143,3421,4213\}$,
$\D^2_{4}=\{2143,3142,4132\}$. We also define
\emph{$\D^1$-Wilf-equivalence} and \emph{$\D^2$-Wilf-equivalence}
similarly to the Wilf-equivalence on $\S_n$. Dumont~\cite{Dumont}
showed that
\[
|\D^1_{2n}|=|\D^2_{2n}|=G_{2n+2}=2(1-2^{2n+2})B_{2n+2},
\]
where $G_n$ is the $n$th Genocchi number, a multiple of the
Bernoulli number $B_n$. Lists of Dumont permutations $\D^1_{2n}$ and
$\D^2_{2n}$ for $n\le 4$ as well as some basic information and
references for Genocchi numbers and Dumont permutations may be
obtained in \cite{Ruskey} and \cite[A001469]{Sloane}. The
exponential generating functions for the unsigned and signed
Genocchi numbers are as follows:
\[
\sum_{n=1}^{\infty}{G_{2n}\frac{x^{2n}}{(2n)!}}=x\tan\frac{x}{2},
\qquad \sum_{n=1}^{\infty}{(-1)^{n}G_{2n}\frac{x^{2n}}{(2n)!}}
=\frac{2x}{e^x+1}-x=-x\tanh\frac{x}{2}.
\]
Some cardinalities of sets of restricted Dumont permutations of
length $2n$ parallel those of restricted permutations of length $n$.
For example, the following results were obtained in
\cite{Bur,Mansour}:
\begin{itemize}
\item $|\D^1_{2n}(\tau)|=C_n$ for $\tau\in\{132,231,312\}$,
where $C_n=\frac{1}{n+1}\binom{2n}{n}$ is the $n$-th Catalan number.

\item $|\D^2_{2n}(321)|=C_n$.
\item $|\D^1_{2n}(213)|=C_{n-1}$, so the operations $r$, $c$ and $r\circ c$ do not
necessarily produce $\D^1$-Wilf-equivalent patterns.
\item $|\D^2_{2n}(231)|=2^{n-1}$, while $|\D^2_{2n}(312)|=1$ and
$|\D^2_{2n}(132)|=|\D^2_{2n}(213)|=0$ for $n\ge 3$, so $r$, $c$ and
$r\circ c$ do not necessarily produce $\D^2$-Wilf-equivalent
patterns either.
\item $|\D^2_{2n}(3142)|=C_n$.
\item
$|\D^1_{2n}(1342,1423)|=|\D^1_{2n}(2341,2413)|=|\D^1_{2n}(1342,2413)|=s_{n+1}$,
the $(n+1)$-st little Schr\"oder number \cite[A001003]{Sloane},
given by $s_1=1$, $s_{n+1}=-s_{n}+2\sum_{k=1}^{n}{s_{k}s_{n-k}}$
($n\ge 2$).
\item $|\D^1_{2n}(2413,3142)|=C(2;n)$, the generalized Catalan
number (see \cite[A064062]{Sloane}).
\end{itemize}

Note that the these results parallel some enumerative avoidance
results in $\S_n$, where the same or similar cardinalities are
obtained:
\begin{itemize}
\item $|\S_n(\tau)|=C_n=\frac{1}{n+1}\binom{2n}{n}$, the $n$th
Catalan number, for any $\tau\in\S_3$.
\item $|\S_n(123,213)|=|\S_n(132,231)|=2^{n-1}$.
\item
$|\S_n(3142,2413)|=|\S_n(4132,4231)|=|\S_n(2431,4231)|=r_{n-1}$, the
$(n-1)$-st large Schr\"oder number \cite[A006318]{Sloane}, given by
$r_0=1$, $r_n=r_{n-1}+\sum_{j=0}^{n-1}{r_k r_{n-k}}$, or
alternatively, by $r_n=2s_{n+1}$.
\end{itemize}

In this paper, we establish several enumerative and bijective
results on restricted Dumont permutations.

In Section \ref{sec:3-letter} we give direct bijections between
$\D_{2n}^1(132)$, $\D_{2n}^1(231)$, $\D_{2n}^2(321)$ and the class
of Dyck paths of semilength $n$ (paths from $(0,0)$ to $(2n,0)$ with
steps $\u=(1,1)$ and $\d=(1,-1)$ that never go below the $x$-axis).
This allows us to consider some permutation statistics, such as
length of the longest increasing (or decreasing) subsequence, and
study their distribution on the sets $\D_{2n}^1(132)$,
$\D_{2n}^1(231)$ and $\D_{2n}^2(321)$.

In Section \ref{sec:4-letter}, we consider Dumont permutations of
the second kind avoiding patterns in $\D^2_4$. Note that \cite{Bur}
showed that $|\D^2_{2n}(3142)|=C_n$ using block decomposition (see
\cite{MV2}), which is very surprising given that it is by far a more
difficult task to count all permutations avoiding a single 4-letter
pattern (e.g., see \cite{Bona,Gessel,Stankova1,Stankova2,West}).

Furthermore, we prove that $\D^2_{2n}(4132)=\D^2_{2n}(321)$ and,
thus, $|\D^2_{2n}(4132)|=C_n$. $\D^2$-Wilf-equivalence of patterns
of different lengths is another striking difference between
restricted Dumont permutations and restricted permutations in
general.

Refining the result $|\D^2_{2n}(3142)|=C_n$ in \cite{Bur}, we
consider some combinatorial statistics on $\D^2_{2n}(3142)$ such as
the number of fixed points and 2-cycles, and give a natural
bijection between permutations in $\D^2_{2n}(3142)$ with $k$ fixed
points and the set $NC(n,n-k)$ of noncrossing partitions of $[n]$
into $n-k$ parts that uses cycle decomposition. This is yet another
surprising difference since pattern avoidance on permutations so far
has not been shown to be related to their cycle decomposition in any
natural way.

Finally, we prove that $|\D^2_{2n}(2143)|=a_na_{n+1}$, where
$a_{2m}=\frac{1}{2m+1}\binom{3m}{m}$ and
$a_{2m+1}=\frac{1}{2m+1}\binom{3m+1}{m+1}$. This allows us to relate
$2143$-avoiding Dumont permutations of the second kind with pairs of
northeast lattice paths from $(0,0)$ to $(2n,n)$ and $(2n+1,n)$ that
do not get above the line $y=x/2$.

Thus, we complete the enumeration problem of $\D^2_{2n}(\tau)$ for
all $\tau\in\D^2_4$. Unfortunately, the same problem for Dumont
permutations of the first kind (i.e. enumeration of permutations in
$\D^1_{2n}(\tau)$ avoiding a pattern in
$\tau\in\D^1_4=\{2143,3421,4213\}$ appears much harder to solve, and
all cases remain unsolved. We do know, however, that no two patterns
in $\D^1_4$ are $\D^1$-Wilf-equivalent \cite{Bur}. On the other
hand, avoidance of pairs of 4-letter patterns yields nice results
\cite{Bur}.

\begin{table}[h]
{\small
\begin{tabular}[h]{ccc}
\begin{tabular}{c|l|l}
  $\tau$ & $|\D^1_{2n}(\tau)|$ & Reference \\
  \hline
  $123$ & Open & \\
  $132$ & $C_n$& \cite[Th. 2.2]{Mansour}\\
  $213$ & $C_{n-1}$& \cite[Th. 2.1]{Bur}\\
  $231$ & $C_n$& \cite[Th. 4.3]{Mansour}\\
  $312$ & $C_n$& \cite[Th. 4.3]{Mansour}\\
  $321$ & $1$& \cite[Page 6]{Bur}\\ \hline
  $2143$& Open & \\
  $3421$& Open & \\
  $4213$& Open & \\ \hline
  $(1342,1423)$& $s_{n+1}$ & \cite[Th. 3.4]{Bur}\\
  $(2341,2413)$& $s_{n+1}$ & \cite[Th. 3.5]{Bur}\\
  $(1342,2413)$& $s_{n+1}$ & \cite[Th. 3.6]{Bur}\\
  $(2341,1423)$& $b_n=3b_{n-1}+2b_{n-2}$ & \cite[Th. 3.7]{Bur}\\
  $(1342,4213)$& $2^{n-1}$ & \cite[Th. 3.9]{Bur}\\
  $(2413,3142)$& $C(2;n)$ & \cite[Th. 3.11]{Bur}\\ \hline
\end{tabular}
& $\qquad$ &
\begin{tabular}{c|l|l}
  $\tau$ & $|\D^2_{2n}(\tau)|$ & Reference\\
  \hline
  $123$ & Open & \\
  $132$ & $0$& Obvious \\
  $213$ & $0$& Obvious \\
  $231$ &$2^{n-1}$& \cite[Th. 2.2]{Bur}\\
  $312$ & $1$& \cite[Page 6]{Bur}\\
  $321$ & $C_n$ & \cite[Th. 4.3]{Mansour}\\ \hline
  $3142$& $C_n$ & \cite[Th. 3.1]{Bur}\\
  $4132$& $C_n$ & Theorem \ref{th4132}\\
  $2143$& $a_na_{n+1}$ & Theorem \ref{th2143}\\ \hline
\end{tabular}
\end{tabular}
}
\medskip
\caption{Some avoidance results for Dumont permutations}
\label{table:known}
\end{table}
Most known avoidance results are given in Table~\ref{table:known}.
Here $a_{2m}=\frac{1}{2m+1}\binom{3m}{m}$ and
$a_{2m+1}=\frac{1}{2m+1}\binom{3m+1}{m+1}$ as defined earlier,
$C(2;n)=\sum_{m=0}^{n-1}{\frac{n-m}{n}\binom{n-1+m}{m}2^m}$, and
$b_0=1$, $b_1=1$, $b_2=3$.

\section{Dumont permutations avoiding a single 3-letter
pattern}\label{sec:3-letter}

In this section we consider some permutation statistics and study
their distribution on certain classes of restricted Dumont
permutations. We focus on the sets $\D_{2n}^1(132)$,
$\D_{2n}^1(231)$ and $\D_{2n}^2(321)$, whose cardinality is given by
the Catalan numbers, as shown in \cite{Bur,Mansour}. We construct
direct bijections between these sets and the class of Dyck paths of
semilength $n$, which we denote $\DD_n$.

\subsection{
$132$-avoiding Dumont permutations of the first kind}
\label{sec:132first}

In this section we present a bijection $\f$ between $\D_{2n}^1(132)$
and $S_n(132)$, which will allow us to enumerate $132$-avoiding
Dumont permutations of the first kind with respect to the length of
the longest increasing subsequences. The bijection is defined as
follows. Let $\pi=\pi_1\pi_2\cdots\pi_{2n}\in\D_{2n}^1(132)$. First
delete all the even entries of $\pi$. Next, replace each of the
remaining entries $\pi_i$ by $(\pi_i+1)/2$. Note that we only obtain
integer numbers since the $\pi_i$ that were not erased are odd.
Clearly, since $\pi$ was $132$-avoiding, the sequence $\f(\pi)$ that
we obtain is a $132$-avoiding permutation, that is,
$\f(\pi)\in\S_n(132)$. For example, if $\pi = 64357821$, then
deleting the even entries we get $3571$, so $\f(\pi)=2341$.

To see that $\f$ is indeed a bijection, we now describe the inverse
map. Let $\sigma\in\S_n(132)$. First replace each entry $\sigma_i$
with $\sigma'_i:=2\sigma_i-1$. Now, for every $i$ from $1$ to $n$,
proceed according to one of the two following cases. If
$\sigma'_i>\sigma'_{i+1}$, insert $\sigma'_i+1$ immediately to the
right of $\sigma'_i$. Otherwise (that is, $\sigma'_i<\sigma'_{i+1}$
or $\sigma'_{i+1}$ is not defined), insert $\sigma'_i+1$ immediately
to the right of the rightmost element to the left of $\sigma'_i$
that is bigger than $\sigma'_i$, or to the beginning of the sequence
if such element does not exist. To see that
$\f^{-1}(\sigma)\in\D_{2n}^1(132)$, note that every even entry
$\sigma'_i+1$ is inserted immediately to the right of either a
smaller odd entry or a larger even entry, or at the beginning of the
sequence, and it is always followed by a smaller entry. Also, after
inserting the even entries, each odd entry $\sigma'_i$ is followed
by an ascent. For example, if $\sigma=546231$, after the first step
we get $(9,7,11,3,5,1)$, so
$\f^{-1}(\sigma)=(9,10,8,7,11,12,4,3,5,6,2,1)$.

Recall Krattenthaler's bijection between $132$-avoiding permutations
and Dyck paths \cite{Krattenthaler}. We denote it
$\kra:\S_n(132)\to\DD_n$, and it can be defined as follows. Given a
permutation $\pi\in\S_n(132)$ represented as an $n\times n$ board,
where for each entry $\pi(i)$ there is a dot in the $i$-th column
from the left and row $\pi(i)$ from the bottom, consider a lattice
path from $(n,0)$ to $(0,n)$ not above the antidiagonal $y=n-x$ that
leaves all dots to the right and stays as close to the antidiagonal
as possible. Then $\kra(\pi)$ is the Dyck path obtained from this
path by reading an $\u$ every time the path goes west and a $\d$
every time it goes north.
 Composing $\f$ with the bijection $\kra$ we obtain a bijection
$\kra\circ\f:\D_{2n}^1(132)\to\DD_n$.

Again through $\kra$, the set $\S_{2n}(132)$ is in bijection with
$\DD_{2n}$. Considering $\D_{2n}^1(132)$ as a subset of
$\S_{2n}(132)$, we observe that
$\inj:=\kra\circ\f^{-1}\circ\kra^{-1}$ is an injective map from
$\DD_n$ to $\DD_{2n}$. Here is a way to describe it directly only in
terms of Dyck paths. Recall that a valley in a Dyck path is an
occurrence of $\d\u$, and that a tunnel is a horizontal segment
whose interior is below the path and whose endpoints are lattice
points belonging to the path (see \cite{Eli,EliPak} for more precise
definitions). Let $D\in\DD_n$. For each valley in $D$, consider the
tunnel whose left endpoint is at the bottom of the valley. Mark the
up-step and the down-step that delimit this tunnel. Now, replace
each unmarked down-step $\d$ with $\d\u\d$. Replace each marked
up-step $\u$ with $\u\u$, and each marked $\d$ with $\d\d$. The path
that we obtain after these operations is precisely
$\inj(D)\in\DD_{2n}$.

To justify the last claim, observe first that a permutation
$\pi\in\D_{2n}^1(132)$ can be decomposed uniquely either as
$\pi=(\tau'+|\tau|,2n-1,2n,\tau)$ or as $\pi=(2n,\tau,2n-1)$, where
$\tau,\tau'$ are again 132-avoiding Dumont permutations of the first
kind, and $|\tau|$ denotes the size of $\tau$. When applying $\kra$
to $\pi\in\D_{2n}^1(132)$, the first decomposition translates into a
Dyck path of the form $C=A\u\u B\d\d$, and the second decomposition
gives a path $C=\u A\d\u\d$, where $A$ and $B$ are Dyck paths. When
the map $\f$ is applied to $\pi$, even entries are deleted, so the
first decomposition becomes
$\f(\pi)=(\f(\tau')+|\f(\tau)|,n,\f(\tau))$, and the second
$\f(\pi)=(\f(\tau),n)$. The translation of this operation in terms
of Dyck paths is that the map $\inj^{-1}=\kra\circ\f\circ\kra^{-1}$
transforms the first decomposition into
$\inj^{-1}(C)=\inj^{-1}(A)\u\inj^{-1}(B)\d$ and the second into
$\inj^{-1}(C)=\u\inj^{-1}(A)\d$. The description of $\inj$ in the
previous paragraph just reverses this construction. Through the map
$\kra$, each entry of the permutation has an associated tunnel in
the path (as described in \cite{Eli}). The construction describing
$\inj$ creates tunnels that correspond to the even elements of
$\f^{-1}(\kra^{-1}(D))$.

For example, if $D=\u\d\u\u\d\u\d\d$, then underlining the marked
steps we get $\u\d\ul{\u}\u\d\ul{\u}\ul{\d}\ul{\d}$, so
$\inj(D)=\u\d\u\d\u\u\u\d\u\d\u\u\d\d\d\d$.

\medskip

Denote by $\lis(\pi)$ (resp. $\lds(\pi)$) the length of the longest
increasing (resp. decreasing) subsequence of $\pi$. Using the above
bijections we obtain the following result.

\begin{theorem}
Let $L_k(z):=\sum_{n\ge0}|\{\pi\in\D_{2n}^1(132):\lis(\pi)\le
k)\}|\,z^n$ be the generating function for
$\{132,12\cdots(k+1)\}$-avoiding Dumont permutations of the first
kind. Then we have the recurrence
\[
L_k(z)=1+\frac{zL_{k-1}(z)}{1-zL_{k-2}(z)},
\]
with $L_{-1}(z)=0$ and $L_0(z)=1$.
\end{theorem}

\begin{proof} As shown in \cite{Krattenthaler}, the length of the longest
increasing subsequence of a permutation $\pi\in\S_{2n}(132)$
corresponds to the height of the path $\kra(\pi)\in\DD_{2n}$. Next
we describe the statistic, which we denote $\lambda$, on the set of
Dyck paths $\DD_n$ that, under the injection
$\inj:\DD_n\hookrightarrow\DD_{2n}$, corresponds to the height in
$\DD_{2n}$. Let $D\in\DD_n$. For each peak $p$ of $D$, define
$\lambda(p)$ to be the height of $p$ plus the number of tunnels
below $p$ whose left endpoint is at a valley of $D$. Now let
$\lambda(D):=\max_p\{\lambda(p)\}$ where $p$ ranges over all the
peaks of $D$. From the description of $\inj$ it follows that for any
$D\in\DD_n$, $\mathrm{height}(\inj(D))=\lambda(D)$. Thus,
enumerating permutations in $\D_{2n}^1(132)$ according to the
parameter $\lis$ is equivalent to enumerating paths in $\DD_n$
according to the parameter $\lambda$. More precisely,
$L_k(z)=\sum_{D\in\DD:\lambda(D)\le k}z^{|D|}$. To find an equation
for $L_k$, we use that every nonempty Dyck path $D$ can be uniquely
decomposed as $D=A\u B\d$, where $A,B\in\DD$. We obtain that
$$L_k(z)=1+zL_{k-1}(z)+z(L_k(z)-1)L_{k-2}(z),$$
where the term $zL_{k-1}(z)$ corresponds to the case where $A$ is
empty (for then $\lambda(\u B\d)=\lambda(B)+1$, and
$z(L_k(z)-1)L_{k-2}(z)$ to the case there $A$ is not empty. From
this we obtain the recurrence
$$L_k(z)=1+\frac{zL_{k-1}(z)}{1-zL_{k-2}(z)},$$
where $L_{-1}(z)=0$ and $L_0(z)=1$ by definition.
\end{proof}

It also follows from the definition of $\kra$ that the length of the
longest decreasing subsequence of $\pi\in\S_{2n}(132)$ corresponds
to the number of peaks of the path $\kra(\pi)\in\DD_{2n}$. Looking
at the description of $\inj$, we see that a peak is created in
$\inj(D)$ for each unmarked down-step of $d$. The number of marked
down-steps is the number of valleys of $D$. Therefore, if
$D\in\DD_n$, we have that the number of peaks of $\inj(D)$ is
$\mathrm{peaks}(\inj(D))=\mathrm{peaks}(D)+n-\mathrm{valleys}(D)=n+1$.
Hence, we have that for every $\pi\in\D_{2n}^1(132)$,
$\lds(\pi)=n+1$.

\subsection{
$231$-avoiding Dumont permutations of the first kind}

As we did in the case of 132-avoiding Dumont permutations, we can
give the following bijection $\ff$ between $\D_{2n}^1(231)$ and
$\S_n(231)$. Let $\pi\in\D_{2n}^1(231)$. First delete all the odd
entries of $\pi$. Next, replace each of the remaining entries
$\pi_i$ by $\pi_i/2$. Note that we only obtain integer entries since
the remaining $\pi_i$ were even. Compare this to the analogous
transformation 
described in Section~\ref{subsec:3142} for Dumont permutations of
the second kind. Clearly the sequence $\ff(\pi)$ that we obtain is a
231-avoiding permutation (since so was $\pi$), that is,
$\ff(\pi)\in\S_n(231)$. For example, if $\pi =
(2,1,10,8,4,3,6,5,7,9)$, then deleting the odd entries we get
$(2,10,8,4,6)$, so $\ff(\pi)=15423$.

To see that $\ff$ is indeed a bijection, we define the inverse map
as follows. Let $\sigma\in\S_n(231)$. First replace each entry $k$
with $2k$. Now, for every $i$ from $1$ to $n-1$, insert $2i-1$
immediately to the left of the first entry to the right of $2i$ that
is bigger than $2i$ (if such an entry does not exist, insert $2i-1$
at the end of the sequence). For example, if $\sigma=7215346$, after
the first step we get $(14,4,2,10,6,8,12)$, so
$\ff^{-1}(\sigma)=(14,4,2,1,3,10,6,5,8,7,9,12,11,13)$.

\medskip

Consider now the bijection $\kra^R:\S_n(231)\longrightarrow\DD_n$
that is obtained by composing $\kra$ defined above with the reversal
operation that sends $\pi=\pi_1\pi_2\cdots\pi_n\in\S_n(231)$ to
$\pi^R=\pi_n\cdots\pi_2\pi_1\in\S_n(132)$.

Through $\kra^R$, the set $\S_{2n}(231)$ is in bijection with
$\DD_{2n}$, so we can identify $\D_{2n}^1(231)$ with a subset of
$\DD_{2n}$. The map $\injj:=\kra^R\circ\ff^{-1}\circ(\kra^R)^{-1}$
is an injection from $\DD_n$ to $\DD_{2n}$. Here is a way to
describe it directly only in terms of Dyck paths. Given $D\in\DD_n$,
all we have to do is replace each down-step $\d$ of $D$ with
$\u\d\d$. The path that we obtain is precisely
$\injj(D)\in\DD_{2n}$. For example, if $D=\u\d\u\u\u\d\u\d\d\d$
(this example corresponds to the same $\sigma$ given above), then
$\injj(D)=\u\u\d\d\u\u\u\u\d\d\u\u\d\d\u\d\d\u\d\d$. Given
$\injj(D)$, one can easily recover $D$ by replacing every $\u\d\d$
by $\d$.

\smallskip

Some properties of $\kra$ trivially translate to properties of
$\kra^R$. In particular, the length of the longest increasing
subsequence of a $231$-avoiding permutation $\pi$ equals the number
of peaks of $\kra^R(\pi)$, and the length of the longest decreasing
subsequence of $\pi$ is precisely the height of $\kra^R(\pi)$.

It follows from the description of $\injj$ in terms of Dyck paths
that for any $D\in\DD_n$, $\injj(D)$ has exactly $n$ peaks (one for
each down-step of $D$). Therefore, for any $\pi\in\D_{2n}^1(231)$,
the number of right-to-left minima of $\pi$ is
$\mathrm{rlm}(\pi)=n$. In fact it is not hard to see directly from
the definition of $231$-avoiding Dumont permutations that the
right-to-left minima of $\pi\in\D_{2n}^1(231)$ are precisely its odd
entries, which necessarily form an increasing subsequence.

Also from the description of $\injj$ we see that
$\mathrm{height}(\injj(D))=\mathrm{height}(D)+1$. In terms of
permutations, this says that if $\pi\in\S_n(231)$, then
$\lds(\ff(\pi))=\lds(\pi)+1$. This allows us to enumerate
$231$-avoiding Dumont permutations with respect to the statistic
$\lds$. Indeed,
$|\{\pi\in\D_{2n}^1(231):\lds(\pi)=k\}|=|\{D\in\DD_n:\mathrm{height}(D)=k-1\}|$.

\subsection{$321$-avoiding Dumont permutations of the second kind}
Let us first notice that a permutation $\pi\in\D_{2n}^2(321)$ cannot
have any fixed points. Indeed, assume that $\pi_i=i$ and let
$\pi=\sigma i \tau$. Since $\pi$ is 321-avoiding, it follows that
$\sigma$ is a permutation of $\{1,2,\ldots,i-1\}$ and $\tau$ is a
permutation of $\{i+1,i+2,\ldots,n\}$. Since $\pi\in\D_{2n}^2$, $i$
must be odd, but then the first element of $\tau$ is in an even
position, and it is either a fixed point or an excedance, which
contradicts the definition of Dumont permutations of the second
kind.

It is known (see e.g.~\cite{Rei}) that a permutation is
$321$-avoiding if and only if both the subsequence determined by its
excedances and the one determined by the remaining elements are
increasing. It follows that a permutation in $\D_{2n}^2(321)$ is
uniquely determined by the values of its excedances. Another
consequence is that if $\pi\in\D_{2n}^2(321)$, then $\lis(\pi)=n$.

We can give a bijection between $\D_{2n}^2(321)$ and $\DD_n$. We
define it in two parts. For the first part, we use the bijection
$\kth$ between $\S_n(321)$ and $\DD_n$ that was defined
in~\cite{Eli}, and which is closely related to the bijection between
$\S_n(123)$ and $\DD_n$ given in~\cite{Krattenthaler}. Given
$\pi\in\S_n(321)$, consider again the $n\times n$ board with a dot
in the $i$-th column from the left and row $\pi(i)$ from the bottom,
for each $i$. Take the path with \emph{north} and \emph{east} steps
that goes from $(0,0)$ to the $(n,n)$, leaving all the dots to the
right, and staying always as close to the diagonal as possible. Then
$\kth(\pi)$ is the Dyck path obtained from this path by reading an
up-step every time the path goes north and a down-step every time it
goes east.

If we apply $\kth$ to a permutation $\pi\in\D_{2n}^2(321)$ we get a
Dyck path $\kth(\pi)\in\DD_{2n}$. The second part of our bijection
is just the map $\injj^{-1}$ defined above, which consists in
replacing every occurrence of $\u\d\d$ with a $\d$. It is not hard
to check that $\pi\mapsto\injj^{-1}(\kth(\pi))$ is a bijection from
$\D_{2n}^2(321)$ to $\DD_n$. For example, for
$\pi=(3,1,5,2,6,4,9,7,10,8)$, we have that
$\kth(\pi)=\u\u\u\d\d\u\u\d\d\u\d\d\u\u\u\d\d\u\d\d$, and
$\injj^{-1}(\kth(\pi))=\u\u\d\u\d\d\u\u\d\d$.

\section{Dumont permutations avoiding a single 4-letter pattern}\label{sec:4-letter}
In this section we will determine the structure
of permutations in $\D_{2n}^2(\tau)$ and find the cardinality
$|\D_{2n}^2(\tau)|$ for each $\tau\in\D^2_4=\{2143,3142,4132\}$.

It was shown in \cite{Bur} that $|\D^2_{2n}(3142)|=C_n$. In Section
\ref{subsec:3142}, we refine this result with respect to the number
of fixed points and 2-cycles in permutations in $\D^2_{2n}(3142)$
and use cycle decomposition to give a natural bijection between
permutations in $\D^2_{2n}(3142)$ with $k$ fixed points and the set
$NC(n,n-k)$ of noncrossing partitions of $[n]$ into $n-k$ parts.

In Section \ref{subsec:4132}, we prove that
$\D^2_{2n}(4132)=\D^2_{2n}(321)$ and, thus, $|\D^2_{2n}(4132)|=C_n$.

Finally, in Section \ref{subsec:2143} we prove that
$|\D^2_{2n}(2143)|=a_na_{n+1}$, where
$a_{2m}=\frac{1}{2m+1}\binom{3m}{m}$ and
$a_{2m+1}=\frac{1}{2m+1}\binom{3m+1}{m+1}$. Thus, we can relate
permutations in $\D^2_{2n}(2143)$ and pairs of northeast lattice
paths from $(0,0)$ to $(n,\left\lfloor\frac{n}{2}\right\rfloor)$ and
$(n+1,\left\lfloor\frac{n+1}{2}\right\rfloor)$ that stay on or below
$y=x/2$.

This completes the enumeration problem of $\D^2_{2n}(\tau)$ for
$\tau\in\D^2_4$.

\subsection{Avoiding 3142}
\label{subsec:3142} 
\begin{figure}[h]
\epsffile{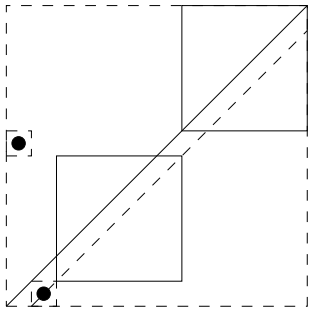} \caption{The block decomposition
of a permutation in $D^2_{2n}(3142)$.} \label{fig:3142-board-2}
\end{figure}

It was shown in \cite{Bur} that
$|\D_{2n}^2(3142)|=C_n$; moreover, the permutations
$\pi\in\D_{2n}^2(3142)$ can be recursively described as follows:
\begin{equation} \label{eq:3142-decomposition}
\pi=(2k,1,r\circ c(\pi')+1,\pi''+2k),
\end{equation}
where $\pi'\in\D_{2k-2}^2(3142)$ and $\pi''\in\D_{2n-2k}^2(3142)$
(see Figure \ref{fig:3142-board-2}). From this block decomposition,
it is easy to see that the subsequence of odd integers in $\pi$ is
increasing. Moreover, the odd entries are exactly those on the main
diagonal and the first subdiagonal (i.e. those $i$ for which
$\pi(i)=i$ or $\pi(i)=i-1$).

In Sections \ref{subsubsec:3142-even} and
\ref{subsubsec:3142-cycles} we use the above decomposition to derive
two bijections from $\D_{2n}^2(3142)$ to sets of cardinality $C_n$.

\medskip

\subsubsection{Subsequence of even entries}
\label{subsubsec:3142-even} The first bijection is
$\phi:\D_{2n}^2(3142)\to E_n\subset S_n$, where
\[
E_n=\left\{(1/2)\pi_{ev}\ |\ \pi\in\D_{2n}^2(3142)\right\},
\]
and $\pi_{ev}$ (resp. $\pi_{ov}$) is the subsequence of even (resp.
odd) values in $\pi$. (Here $\frac{1}{2}\pi_{ev}$ is the permutation
obtained by dividing all entries in $\pi_{ev}$ by 2; in other words,
if $\sigma=\frac{1}{2}\pi_{ev}$, then $\sigma(i)=\pi_{ev}(i)/2$ for
all $i\in[n]$.) Define $\phi(\pi)=\frac{1}{2}\pi_{ev}$ for each
$\pi\in\D_{2n}^2(3142)$.

Permutations in $E_n$ have a block decomposition similar to those in
$\D_{2n}^2(3142)$, namely,
\[
\sigma\in E_n \iff \sigma=(k,r\circ c(\sigma'),k+\sigma'')\
\text{for some } \sigma'\in E_{k-1}\ \text{and}\ \sigma''\in
E_{n-k}.
\]

The inverse $\phi^{-1}:E_n\to\D^2_{2n}(3142)$ is easy to describe.
Let $\sigma\in E_n$. Then $\pi=\phi^{-1}(\sigma)$ is obtained as
follows: let $\pi_{ev}=2\sigma$ (i.e. $\pi_{ev}(i)=2\sigma(i)$ for
all $i\in[n]$), then for each $i\in[n]$ insert $2i-1$ immediately
before $2\sigma(i)$ if $\sigma(i)<i$ or immediately after
$2\sigma(i)$ if $\sigma(i)\ge i$. For instance, if $\sigma=3124\in
E_4$, then $\pi_{ev}=6248$ and
$\pi=61\,32\,54\,87\in\D^2_{8}(3142)$.

It is not difficult to show that $E_n$ consists of exactly those
permutations that, written in cyclic form, correspond to noncrossing
partitions of $[n]$ by replacing pairs of parentheses with slashes.
We remark that $E_n$ is also the set of permutations whose tableaux
(see \cite{SW}) have a single 1 in each column.

\begin{theorem} \label{thm:3142_e}
For a permutation $\rho$, define
\begin{xalignat*}{3}
\mathrm{fix}(\rho)&=|\{i\ |\ \rho(i)=i\}|, &\mathrm{exc}(\rho)&=|\{i\ |\ \rho(i)>i\}|,\\
\mathrm{fix}_{_{-1}}(\rho)&=|\{i\ |\ \rho(i)=i-1\}|,
&\mathrm{def}(\rho)&=|\{i\ |\ \rho(i)<i\}|.
\end{xalignat*}
Then for any $\pi\in\D_{2n}^2(3142)$ and $\sigma=\phi(\pi)\in E_n$,
we have
\begin{align}
\mathrm{fix}(\pi)+\mathrm{fix}_{_{-1}}(\pi)&=n, \label{eq:3142-n}\\
\mathrm{fix}(\pi)&=\mathrm{def}(\sigma), \label{eq:3142-fix}\\
\mathrm{fix}_{_{-1}}(\pi)&=\mathrm{exc}(\sigma)+\mathrm{fix}(\sigma), \label{eq:3142-fix-1}\\
\mathrm{fix}(\sigma)&=\#\ \text{\rm 2-cycles in }\pi.
\label{eq:3142-fix-2}
\end{align}
\end{theorem}

\begin{proof}
Equation (\ref{eq:3142-n}) follows from the fact that odd integers
in $\pi$ are exactly those on the main diagonal and first
subdiagonal.

Let $\pi$ and $\sigma$ be as above and let $i\in[n]$. Then there are
two cases: either $2i-1=\pi(2i)$ or $2i-1=\pi(2i-1)$.

\emph{Case 1:} $\pi(2i)=2i-1$. Then $\pi(2i-1)\ge 2i$, and hence
$\pi(2i-1)$ must be even.

\emph{Case 2:} $\pi(2i-1)=2i-1$. Then $\pi(2i)\le 2i-2$, and hence
$\pi(2i)$ must be even.\\ In either case, for each $i\in [n]$, we
have $\{\pi(2i-1),\pi(2i)\}=\{2i-1,2s_i\}$ for some $s_i\in [n]$.
Define $\sigma(i)=s_i$. Then $\sigma(i)\ge i$ if $2i-1\in
\mathrm{fix}_{_{-1}}(\pi)$, and $\sigma(i)\le i-1$ if $2i-1\in
\mathrm{fix}(\pi)$. This proves (\ref{eq:3142-fix}) and
(\ref{eq:3142-fix-1}).

Finally, let $i\in[n]$ be such that $\sigma(i)=i$. Since
$2\sigma(i)\in\{\pi(2i-1),\pi(2i)\}$ and $\pi(2i)<2i$, it follows
that $2i=2\sigma(i)=\pi(2i-1)$, so $2i-1=\pi(2i)$, and thus $\pi$
contains a 2-cycle $(2i-1,2i)$.

Conversely, let $(ab)$ be a 2-cycle of $\pi$, and assume that $b>a$.
Then $\pi(a)>a$, so $a$ must be odd, say $a=2i-1$ for some
$i\in[n]$. Then $b=\pi^{-1}(a)\in\{2i-1,2i\}$, so $b=2i$, and thus
$(ab)=(2i-1,2i)$. This proves (\ref{eq:3142-fix-2}).
\end{proof}

\begin{theorem} \label{thm:3142-fix-2}
Let $A(q,t,x)=\sum_{n\ge0}\sum_{\pi\in \D^2_{2n}(3142)}
q^{\mathrm{fix}(\pi)} t^{\# \text{ 2-cycles in }\pi}x^n$ be the
generating function for $3142$-avoiding Dumont permutations of the
second kind with respect to the number of fixed points and the
number of 2-cycles. Then
\begin{equation}\label{eq:A3142}
A(q,t,x)=\frac{1+x(q-t)-\sqrt{1-2x(q+t)+x^2((q+t)^2-4q)}}{2xq(1+x(1-t))}.
\end{equation}
\end{theorem}
\begin{proof}
By the correspondences in Theorem~\ref{thm:3142_e}, it follows that
$$A(q,t,x)=\sum_{n\ge0}\sum_{\sigma\in E_n} q^{\mathrm{def}(\sigma)}
t^{\mathrm{fix}(\sigma)}x^n.$$ For convenience, let us define a
related generating function $B(q,t,x)=\sum_{n\ge0}\sum_{\sigma\in
E_n} q^{\mathrm{def}(\sigma)} t^{\mathrm{fix}_{_{-1}}(\sigma)}x^n$.
From the block decomposition of permutations $\sigma\in E_n$ as
$\sigma=(k,r\circ c(\sigma'),k+\sigma'')$ for some $\sigma'\in
E_{k-1}$, $\sigma''\in E_{n-k}$, it follows that
\begin{equation}
A(q,t,x)=1+xt A(q,t,x) + x (B(1/q,t,xq)-1) A(q,t,x). \label{eq:A}
\end{equation}
The term $xt A(q,t,x)$ corresponds to the case $k=1$, in which
$\sigma'$ is empty and $k$ is a fixed point. When $k>1$, $\sigma''$
still contributes as $A(q,t,x)$, and the contribution of $\sigma'$
is $B(1/q,t,xq)-1$, since elements with $\sigma'(i)=i-1$ become
fixed points of $\sigma$, and all elements of $\sigma'$ other than
its deficiencies become deficiencies of $\sigma$.

A similar reasoning gives the following equation for $B(q,t,x)$:
$$B(q,t,x)=1+x A(1/q,t,xq) B(q,t,x).$$
Solving for $B$ we have $B(q,t,x)=\frac{1}{1-xA(1/q,t,xq)}$, and
plugging $B(1/q,t,xq)=\frac{1}{1-xqA(q,t,x)}$ into (\ref{eq:A})
gives
$$A(q,t,x)=1+x\left(\frac{1}{1-xqA(q,t,x)}+t-1\right)A(q,t,x).$$
Solving this quadratic equation gives the desired formula for
$A(q,t,x)$.
\end{proof}

\medskip

\subsubsection{Cycle decomposition} \label{subsubsec:3142-cycles} Letting $t=1$
in (\ref{eq:A3142}), we obtain
\begin{corollary}\label{cor:3142-noncrossing}
We have
\[
\sum_{n\ge0}\sum_{\pi\in \D^2_{2n}(3142)} q^{\mathrm{fix}(\pi)} x^n=
A(q,1,x)=\frac{1+x(q-1)-\sqrt{1-2x(q+1)+x^2(q-1)^2}}{2xq},
\]
i.e. the number of permutations in $\pi\in\D^2_{2n}(3142)$ with $k$
fixed points is the Narayana number
$N(n,k)=\frac{1}{n}\binom{n}{k}\binom{n}{k+1}$, which is also the
number of noncrossing partitions of $[n]$ into $n-k$ parts.
\end{corollary}

\begin{figure}[h]
\epsfxsize=160.0pt \epsffile{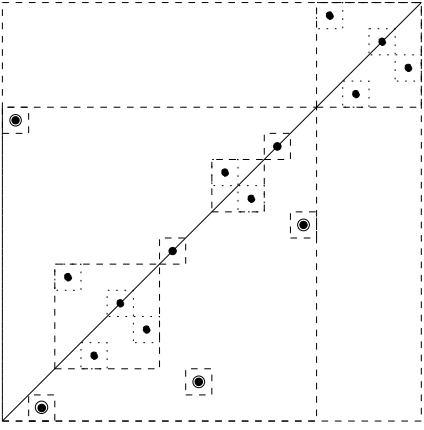} \caption{The
cycle decomposition of a permutation in $D^2_{2n}(3142)$. The
circled dots correspond to one of the cycles.}
\label{fig:3142-cycle}
\end{figure}

\begin{proof}
Even though the generating function is an immediate consequence of
Theorem~\ref{thm:3142-fix-2}, we will give a combinatorial proof of
the corollary, by exhibiting a natural bijection
$\psi:\D_{2n}^2(3142)\to NC(n)$, where $NC(n)$ is the set of
noncrossing partitions of $[n]$. We start by considering a
permutation $\pi\in\D_{2k}^2(3142)$. Iterating the block
decomposition (\ref{eq:3142-decomposition}), we obtain
\[
\begin{split}
\pi&=(2k_1, 1, c\circ r(\pi_1)+1, 2k_2, 2k_1+1, c\circ
r(\pi_2)+2k_1+1, \cdots, 2k_r, 2k_{r-1}+1, c\circ
r(\pi_r)+2k_{r-1}+1)\\
&=(2k_1, 1, 2k_1-r(\pi_1), 2k_2, 2k_1+1, 2k_2-r(\pi_2), \cdots,
2k_r, 2k_{r-1}+1, 2k_{r}-r(\pi_r)),
\end{split}
\]
where $1\le k_1<k_2<\dots<k_r=k$,
$\pi_i\in\D_{2(k_i-k_{i-1}-1)}^2(3142)$ ($1\le i\le r$), and we
define $k_0=0$. Note that each permutation $c\circ
r(\pi_i)+2k_{i-1}+1=2k_i-r(\pi_i)$ of $[2k_{i-1}+2,2k_i-1]$ occurs
at positions $[2k_{i-1}+3,2k_i]$ in $\pi$.

Now consider
\[
\pi'=(2k+2,1,c\circ r(\pi)+1)=(2k_r+2,1,2k_r+2-r(\pi)).
\]
Let $k'_i=k-k_i=k_r-k_i$. By (\ref{eq:3142-decomposition}), we have
$\pi\in\D_{2k+2}^2(3142)$, $\pi_i\in
\D_{2(k'_{i-1}-k'_i-1)}^2(3142)$ ($1\le i\le r$), $k'_r=0$,
$k'_0=k$, and
\[
\pi'=(2k+2,1,\pi_r+2,2k'_{r-1}+1,2,\pi_{r-1}+2k'_{r-1}+2,
2k'_{r-2}+1,2k'_{r-1}+2, \dots, \pi_1+2k'_1+2,2k+1,2k'_1+2).
\]

Note that, for each $i=1,2,\dots,r$, the permutation $\pi_i+2k'_i+2$
of $[2k'_i+3,2k'_{i-1}]$ occurs at positions $[2k'_i+3,2k'_{i-1}]$
in $\pi'$. Moreover, the entries $2k'_i+1$ ($0\le i\le r-1$) occur
at positions $2k'_i+1$ in $\pi'$, and thus are fixed points of
$\pi'$. Finally, each entry $2k'_i+2$ ($1\le i\le r$) occurs at
position $2k'_{i-1}+2$, 1 occurs at position $2=2k'_r+2$, and
$2k+2=2k'_0+2$ occurs at position $1$.

Thus, $\gamma=(2k'_0+2,2k'_1+2,2k'_2+2,\dots,2k'_{r-1}+2,2k'_r+2,1)
=(2k+2,2k'_1+2,2k'_2+2,\dots,2k'_{r-1}+2,2,1)$ is a cycle of $\pi'$
(such as the one consisting of circled dots in
Figure~\ref{fig:3142-cycle}), and each remaining nontrivial cycle of
$\pi'$ is completely contained in some $\pi_i+2k'_i+2$, which is a
$3142$-avoiding Dumont permutation of the second kind of
$[2k'_i+3,2k'_{i-1}]$. Note that
\[
2k'_i+2<2k'_i+3<2k'_{i-1}<2k'_{i-1}+2,
\]
so all entries of each remaining cycle of $\pi'$ are contained
between two consecutive entries of $\gamma$.

Now let $G$ be the subset of $[2k+2]$ consisting of the entries of
$\gamma$. Then, clearly,
\[
G/\{2k'_{r-1}+1\}/\dots/\{2k'_1+1\}/\{2k'_0+1\}/[2k'_r+3,2k'_{r-1}]/\dots/[2k'_1+3,2k'_{0}]
\]
is a noncrossing partition of $[2k+2]$. Now it is easy to see by
induction on the size of $\pi'$ that the subsets of $\pi'$ formed by
entries of the cycles in cycle decomposition of $\pi'$ form a
noncrossing partition of $\pi'$. Moreover, all the entries of $G$
except the smallest entry are even, so likewise the cycle
decomposition of $\pi'$ determines a unique noncrossing partition of
$\pi'_{ev}$, hence a unique noncrossing partition of $[n]$.

Finally, any permutation $\hat\pi\in\D_{2n}^2(3142)$ can be written
as $\hat\pi=(\pi',\pi''+2k+2)$, where $\pi'$ is as above and
$\pi''\in\D_{2n-2k-2}^2(3142)$, so the cycles of any permutation in
$\D_{2n}^2(3142)$ determine a unique noncrossing partition of $[n]$.

Notice also that each cycle in the decomposition of $\hat\pi$
contains exactly one odd entry, the least entry in each cycle, so
the number of odd entries of $\hat\pi$ which are not fixed points,
$\mathrm{fix}_{_{-1}}(\hat\pi)=n-\mathrm{fix}(\hat\pi)$, is the
number of parts in $\psi(\hat\pi)$. This finishes the proof.
\end{proof}

For example (see Figure~\ref{fig:3142-cycle}), if
\[
\begin{split}
\hat\pi&=12,1,6,3,5,4,7,2,10,9,11,8,16,13,15,14\\
&=(12,8,2,1)(6,4,3)(10,9)(16,14,13)(15)(11)(7)(5)\in\D_{16}^2(3142),
\end{split}
\]
then $\psi(\hat\pi)=641/32/5/87\in NC(8)$. Note also that
$\hat\pi_{ev}=63215487=(641)(32)(5)(87)$.

\subsection{Avoiding 4132}
\label{subsec:4132} For Dumont permutations of the second kind
avoiding the pattern $4132$ we have the following result.

\begin{theorem}\label{th4132}
For any $n\ge0$, $\D_{2n}^2(4132)=\D_{2n}^2(321)$. Moreover,
$|\D_{2n}^2(4132)|=C_n$, where $C_n$ is the $n$th Catalan number.
Thus, $4132$ and $3142$ are $\D^2$-Wilf-equivalent.
\end{theorem}
\begin{proof}
The pattern $321$ is contained in $4132$. Therefore, if $\pi$ avoids
$321$, then $\pi$ avoids $4132$, so $\D_{2n}^2(321)\subseteq
\D_{2n}^2(4132)$. Now let us prove that $\D_{2n}^2(4132)\subseteq
\D_{2n}^2(321)$. Let $n\ge4$ and let $\pi\in \D_{2n}^2(4132)$
contain an occurrence of $321$. Choose the occurrence of $321$ in
$\pi$, say $\pi(i_1)>\pi(i_2)>\pi(i_3)$ with $1\le i_1<i_2< i_3\le
2n$, such that $i_1+i_2+i_3$ is minimal. If $i_1$ is an even number,
then $\pi(i_1-1)\ge i_1-1\ge \pi(i_1)$, so the occurrence
$\pi(i_1-1)\pi(i_1)\pi(i_2)$ of pattern $321$ contradicts minimality
of our choice. Therefore, $i_1$ is odd. If $i_2\ne i_1+1$, then from
the minimality of the occurrence we get that $\pi(i_1+1)<\pi(i_3)$.
Hence, $\pi$ contains $4132$, a contradiction. So $i_2=i_1+1$. If
$i_3$ is odd, then $\pi(i_3)\ge i_3>i_1+1\ge \pi(i_1+1)$, which
contradicts $\pi(i_1)>\pi(i_1+1)>\pi(i_3)$. So $i_3$ is even.

Therefore, the our chosen occurrence of $321$ is given by
$\pi(2i+1)\pi(2i+2)\pi(j)$ where $4\le 2i+2\le j\le 2n$ (since
$\pi(2)=1$, we must have $i\ge1$). By minimality of the occurrence,
we have $\pi(m)\le 2i$ for all $m\le 2i$. On the other hand,
$\pi(i_3)<\pi(2i+2)\le 2i+1$ which means that $\pi(i_3)\le 2i$.
Hence, $\pi$ must contain at least $2i+1$ letters smaller than $2i$,
a contradiction.

Thus, if $\pi\in \D_{2n}^2(4132)$ then $\pi\in \D_{2n}^2(321)$. The
rest follows from \cite[Theorem 4.3]{Mansour}.
\end{proof}

\subsection{Avoiding 2143}
\label{subsec:2143} Dumont permutations of the second kind that
avoid $2143$ are enumerated by the following theorem, which we prove
in this section.

\begin{theorem}\label{th2143}
For any $n\ge0$, $|\D_{2n}^2(2143)|=a_{n}a_{n+1}$, where
\[
\begin{split}
a_{2m}&=\frac{1}{2m+1}\binom{3m}{m},\\
a_{2m+1}&=\frac{1}{2m+1}\binom{3m+1}{m+1}=\frac{1}{m+1}\binom{3m+1}{m}.
\end{split}
\]
\end{theorem}

\begin{remark} \label{rem:ternary}
Note that the sequence $\{a_n\}$ also enumerates northeast lattice
paths in $\mathbb{Z}^2$ from $(0,0)$ to
$(n,\left\lfloor\frac{n}{2}\right\rfloor)$ that stay on or below
$y=x/2$, as well as symmetric ternary trees on $3n$ edges and
symmetric diagonally convex directed polyominoes with $n$ squares
(see \cite[A047749]{Sloane} and references therein). Also note that
$\{a_{2m+1}\}$ is the convolution of $\{a_{2m}\}$ with itself, while
the convolution of $\{a_{2m}\}$ with $\{a_{2m+1}\}$ is
$\{a_{2m+2}\}$. Alternatively, if $f(x)$ and $g(x)$ are the ordinary
generating functions for $\{a_{2m}\}$ and $\{a_{2m+1}\}$, then
$f(x)=1+xf(x)g(x)$ and $g(x)=f(x)^2$, so $f(x)=1+xf(x)^3$. Now the
Lagrange inversion applied to the last two equations yields the
formulas for $a_n$.
\end{remark}

\begin{lemma} \label{lem:odd-even-split}
Let $\pi\in\D^2_{2n}(2143)$. Then the subsequence
$(\pi(1),\pi(3),\dots,\pi(2n-1))$ is a permutation of
$\{n+1,n+2,\dots,2n\}$ and the subsequence
$(\pi(2),\pi(4),\dots,\pi(2n))$ is a permutation of
$\{1,2,\dots,n\}$.
\end{lemma}

\begin{proof}
Assume the lemma is false. Let $i$ be the smallest integer such that
$\pi(2i)\ge n+1$. Then $\pi(2i-1)\ge 2i-1\ge\pi(2i)\ge n+1$.
Therefore, if $j\ge i$, then $\pi(2j-1)\ge 2j-1\ge 2i-1\ge n+1$. In
fact, note that for any $1\le j\le n$, $\pi(2j-1)\ge 2j-1\ge
\pi(2j)$.

By minimality of $i$, we have $\pi(2j)\le n$ for $j<i$. Hence, if
$\pi(2j-1)\le n$ for some $j<i$, then
$(\pi(2j-1),\pi(2j),\pi(2i-1),\pi(2i))$ is an occurrence of pattern
$2143$ in $\pi$. Therefore, $\pi(2j-1)\ge n+1$ for all $j<i$.

Thus, we have $\pi(2j-1)\ge n+1$ for any $1\le j\le n$, and
$\pi(2i)\ge n+1$, so $\pi$ must have at least $n+1$ entries between
$n+1$ and $2n$, which is impossible. The lemma follows.
\end{proof}

For $\pi\in\D^2_{2n}(2143)$, we denote
$\pi_o=(\pi(1),\pi(3),\dots,\pi(2n-1))-n$ and
$\pi_e=(\pi(2),\pi(4),\dots,\pi(2n))$. By Lemma
\ref{lem:odd-even-split}, $\pi_o,\pi_e\in\S_n(2143)$. For example,
given $\pi=71635482\in\D^2_8(2143)$, we have $\pi_o=3214$ and
$\pi_e=1342$. Note that $\pi(2i-1)=\pi_o(i)+n$ and
$\pi(2i)=\pi_e(i)$.

\begin{figure}[h]
\epsffile{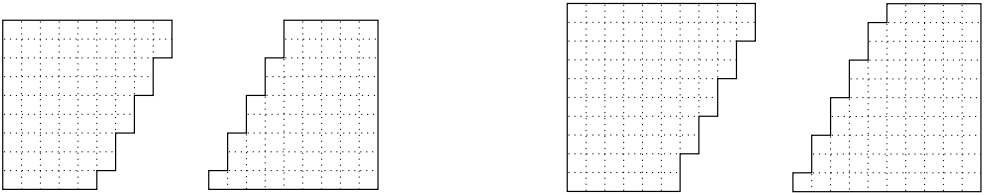} \caption{The boards of Lemma
\ref{lem:boards} for $n=9$ (left) and $n=10$ (right).}
\label{fig:2143-boards}
\end{figure}

\begin{lemma}\label{lem:boards}
For any permutation $\pi\in\D^2_{2n}(2143)$, and $\pi_o$ and $\pi_e$
defined as above, the following is true:\\[-16pt]
\begin{enumerate}
\item\label{lem:boards-odd} $\pi_o\in\S_n(132)$ and the entries of
$\pi_o$ are on a board with $n$ top-justified columns of sizes
$2,4,6,\dots,2\left\lfloor\frac{n}{2}\right\rfloor,n,\dots,n$ from
right to left (see the first and third boards in Figure
\ref{fig:2143-boards}).

\item\label{lem:boards-even} $\pi_e\in\S_n(213)$ and the entries
of $\pi_e$ are on a board with $n$ bottom-justified columns of sizes
$1,3,5,\dots,2\left\lfloor\frac{n}{2}\right\rfloor-1,n,\dots,n$ from
left to right (see the second and fourth boards in Figure
\ref{fig:2143-boards}).

\end{enumerate}
\end{lemma}

\begin{proof}
If $132$ occurs in $\pi_o$ at positions $i_1<i_2<i_3$, then $2143$
occurs in $\pi$ at positions $2i_1-1<2i_1<2i_2-1<2i_3-1$ since
$\pi(2i_1)<\pi(2i_1-1)$. Similarly, if $213$ occurs in $\pi_e$ at
positions $i_1<i_2<i_3$, then $2143$ occurs in $\pi$ at positions
$2i_1<2i_2<2i_3-1<2i_3$ since $\pi(2i_3-1)>\pi(2i_3)$. The rest
simply follows from the definition of $\D^2_{2n}$ and Lemma
\ref{lem:odd-even-split}.
\end{proof}

Let us call a permutation as in part (\ref{lem:boards-odd}) of Lemma
\ref{lem:boards} an \emph{upper board}, and a permutation as in part
(\ref{lem:boards-even}) of Lemma \ref{lem:boards} a \emph{lower
board}. Note that $\pi_e(1)=1$ and $213=r\circ c(132)$. Hence it is
easy to see that $\pi_e=(1,r\circ c(\pi')+1)$ with
$\pi'\in\S_{n-1}(132)$ of upper type. Let $b_n$ be the number of
lower boards in $\S_n(213)$. Then the number of upper boards in
$\S_n(132)$ is $b_{n+1}$.

\begin{lemma} \label{lem:2143-u-l}
Let $\pi_1\in\S_n(132)$ be an upper board and $\pi_2\in\S_n(213)$ be
a lower board. Let $\pi\in\S_{2n}$ be defined by
$\pi=(\pi_1(1)+n,\pi_2(1),
\pi_1(2)+n,\pi_2(2),\dots,\pi_1(n)+n,\pi_2(n))$ (i.e. such that
$\pi_o=\pi_1$ and $\pi_e=\pi_2$). Then $\pi\in\D^2_{2n}(2143)$.
\end{lemma}

\begin{figure}[h]
\epsfxsize=120.0pt \epsffile{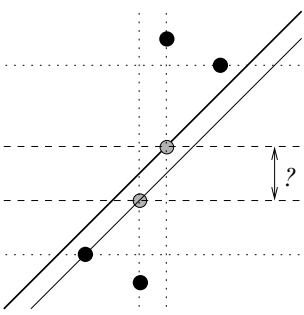} \caption{This
situation is impossible in Lemma \ref{lem:2143-u-l}: no value
between the grey points (inclusive) can occur in $\pi$.}
\label{fig:21-43}
\end{figure}

\begin{proof}
Clearly $\pi\in\D^2_{2n}$. It is not difficult to see that if $\pi$
contains $2143$, then ``2'' and ``1'' are deficiencies (i.e., they
are at even positions and come from $\pi_2$) and ``4'' and ``3'' are
excedances or fixed points (i.e. they are at odd positions and come
from $\pi_1$). Such an occurrence is represented in Figure
\ref{fig:21-43}, where an entry $\pi(i)$ is plotted by a dot with
abscissa $i$ and ordinate $\pi(i)$, and the two diagonal lines
indicate the positions of the fixed points and elements with
$\pi(i)=i-1$.

Say the pattern $2143$ occurs at positions
$2i_1<2i_2<2i_3-1<2i_4-1$. We have $\pi(2j)\le 2j-1< 2i_2-1$ for any
$j<i_2$. On the other hand, the subdiagonal part of $\pi$ avoids
$213$, so $\pi(2j)<\pi(2i_1)\le 2i_1-1<2i_2-1$ for any $j\ge i_2$.
Thus, $\pi(2j)<2i_2-1$ for any $1\le j\le n$. Similarly,
$\pi(2j-1)\ge 2j-1>2i_3-1$ for any $j>i_3$, and
$\pi(2j-1)>\pi(2i_4)\ge 2i_4-1>2i_3-1$ for any $j\le i_3$ since the
superdiagonal part of $\pi$ avoids $132$. Thus, $\pi(2j-1)>2i_3-1$
for any $1\le j\le n$.

Therefore, no entry of $\pi$ lies in the interval $[2i_2-1,2i_3-1]$,
which is nonempty since $2i_2<2i_3-1$. This is, of course,
impossible, so the lemma follows.
\end{proof}

Hence, there is a 1-1 correspondence between permutations
$\pi\in\D^2_{2n}(2143)$ and pairs of permutations $(\pi_1,\pi_2)$,
where $\pi_1\in\S_n(132)$ is an upper board and $\pi_2\in\S_n(213)$
is a lower board. Thus, $|\D^2_{2n}(2143)|=b_{n}b_{n+1}$, where
$b_n$ is the number of lower boards $\pi\in\S_n(213)$ and $b_{n+1}$
is the number of upper boards $\pi\in\S_n(132)$ (see the paragraph
before Lemma \ref{lem:2143-u-l}).

\begin{lemma} \label{lem:2143-lower}
Let $F(x)=\sum_{m=0}^{\infty}{b_{2m}x^m}$ and
$G(x)=\sum_{m=0}^{\infty}{b_{2m+1}x^{m}}$. Then we have $b_0=1$ and
\[
\begin{split}
b_{2m}=\sum_{i=0}^{m-1}{b_{2i}b_{2m-2i-1}}, &\quad
b_{2m+1}=\sum_{i=0}^{m}{b_{2i}b_{2m-2i}},\\
F(x)=1+xF(x)G(x), &\quad G(x)=F(x)^2.
\end{split}
\]
\end{lemma}

\begin{figure}[h]
\epsffile{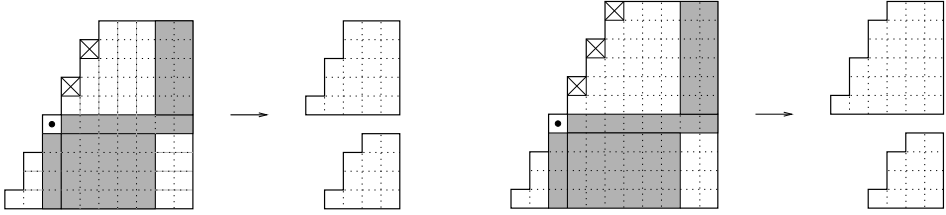} \caption{A lower board
$\pi\in\S_n(213)$ ($n=10$ (even), left, and $n=11$ (odd), right)
decomposed into two lower boards according to the largest $i$ such
that $\pi(i+1)=2i+1$ (here $i=2$). } \label{fig:2143-lower}
\end{figure}

\begin{proof}
Let $\pi\in S_n(213)$ be a lower board, and let $i\ge 0$ be maximal
such that $\pi(i+1)=2i+1$. Such an $i$ always exists since
$\pi(1)=1$. Then $\pi(j)\le 2j-2$ for $j\ge i+2$. Furthermore, $\pi$
avoids $213$, so if $j_1,j_2>i+1$, and $\pi(j_1)>\pi(i+1)>\pi(j_2)$,
then $j_1<j_2$. In other words, all entries of $\pi$ greater than
and to the right of $2i+1$ must come before all entries less than
and to the right of $2i+1$ (see Figure \ref{fig:2143-lower}, the
areas that cannot contain entries of $\pi$ are shaded). In addition,
$\pi(j)\le 2i+1$ for $j\le i+1$, so $\pi(j)>2i+1$ only if $j>i+1$.
There are $n-2i-1$ values greater than $2i+1$ in $\pi$, hence they
must occupy the $n-2i-1$ positions immediately to the right of
$\pi(i+1)$, i.e. positions $i+2$ through $n-i$. It is not difficult
now to see from the above argument that all entries of $\pi$ greater
than $2i+1$ must lie on a board of lower type in $\S_{n-2i-1}(213)$,
while the entries less than $2i+1$ in $\pi$ must lie on two boards
whose concatenation is a lower board in $\S_{2i}(213)$ (unshaded
areas in Figure \ref{fig:2143-lower}).
\end{proof}

Thus, we get the same generating function equations as in Remark
\ref{rem:ternary}, so $F(x)=f(x)$, $G(x)=g(x)$, and hence $b_n=a_n$
for all $n\ge 0$. This proves Theorem \ref{th2143}.

We can give a direct bijection showing that $b_n=a_n$. It is
well-known that $a_{2n}$ (resp. $a_{2n+1}$) is the number of
northeast lattice paths from $(0,0)$ to $(2n,n)$ (resp. from $(0,0)$
to $(2n+1,n)$) that do not get above the line $y=x/2$. The following
bijection uses the same idea as a bijection of Krattenthaler
\cite{Krattenthaler} from the set of $132$-avoiding permutations in
$\S_n$ to Dyck paths of semilength $n$, which is described in
Section~\ref{sec:132first}.

\begin{figure}[h]
\epsffile{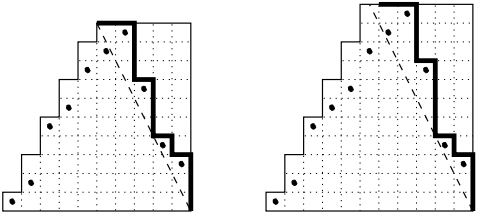} \caption{A bijection between lower
boards in $\S_n(213)$, for $n=10$ (left) and $n=11$ (right), and
northwest paths from $(n,0)$ to $(\lceil n/2 \rceil, n)$ not below
$y=2n-2x$.} \label{fig:2143_path}
\end{figure}

We introduce a bijection between the set of lower boards in
$\S_n(213)$ and northwest paths from $(n,0)$ to $(\lceil n/2 \rceil,
n)$ that stay on or above the line $y=2n-2x$ (see Figure
\ref{fig:2143_path}). Given a lower board in $\S_n(213)$ represented
as an $n\times n$ binary array, consider a lattice path from $(n,0)$
to $(\lceil n/2 \rceil, n)$ that leaves all dots to the left and
stays as close to the $y=2n-2x$ as possible. We claim that such a
path must stay on or above the line $y=2n-2x$. Indeed, considering
rows of a lower board from top to bottom, we see that at most one
extra column appears on the left for every two consecutive rows.
Therefore, our path must shift at least $r$ columns to the right for
every $2r$ consecutive rows starting from the top. The rest is easy
to see.

Conversely, given a northwest path from $(n,0)$ to $(\left\lceil n/2
\right\rceil, n)$ not below the line $y=2n-2x$, fill the
corresponding board from top to bottom (i.e. from row $n$ to row
$1$) so that the dots are in the rightmost column to the left of the
path that still contains no dots.

Theorem \ref{th2143} implies that
$\lim_{n\to\infty}{|\D_{2n}^2(2143)|^{\frac{1}{2n}}=\frac{3^3}{2^2}=\frac{27}{4}}$.
In comparison, \cite{Regev} and \cite{West} imply that
$|\S_{n}(2143)|=|\S_{n}(1234)|$ and hence
$\lim_{n\to\infty}{|\S_{n}(2143)|^\frac{1}{n}}=\lim_{n\to\infty}{|\S_{n}(1234)|^\frac{1}{n}}=(4-1)^2=9$.

The \emph{median Genocchi number} (or \emph{Genocchi number of the
second kind}) $H_n$ \cite[A005439]{Sloane} counts the number of
derangements in $\D^2_{2n}$ (also, the number of permutations in
$\D^1_{2n}$ which begin with $n$ or $n+1$). Using the preceding
argument, we can also count the number of derangements in
$\D^2_{2n}(2143)$.

\begin{theorem} \label{th2143median}
The number of derangements in $\D^2_{2n}(2143)$ is $a_n^2$, where
$a_n$ is as in Theorem \ref{th2143}.
\end{theorem}

\begin{proof}
Notice that the fixed points of a permutation
$\pi\in\D^2_{2n}(2143)$ correspond to the dots in the lower right
(southeast) corner cells on its upper board (except the lowest right
corner when $n$ is odd) (see Figure \ref{fig:2143-boards}). It is
easy to see that deletion of those cells on an upper board produces
a rotation of a lower board by $180^\circ$. This, together with the
preceding lemmas, implies the theorem.
\end{proof}

The following theorem gives the generating function for the
distribution of the number of fixed points among permutations in
$\D^2_{2n}(2143)$.
\begin{theorem} \label{thm:2143-fix}
We have
\begin{equation} \label{eq:2143-fix}
\begin{split}
\sum_{\pi\in\D^2_{2n}(2143)}{q^{\mathrm{fix}(\pi)}}&=a_n\cdot
[x^{n+1}]\left(\frac{1}{1-xf(x^2)}\cdot\frac{1}{1-qx^2f(x^2)^2}\right)\\
&=a_n\cdot[x^{n+1}]\frac{f(x^2)}{(1-xf(x^2))(q+(1-q)f(x^2))}.
\end{split}
\end{equation}
where $f(x)=\sum_{n\ge 0}{a_{2n}x^n}$ is a solution of
$f(x)=1+xf(x)^3$, and $[x^n]h(x)$ is the coefficient at $x^n$ in the
power series representation of $h(x)$.
\end{theorem}

Note that $\sum_{n\ge 0}{a_{2n}x^{2n}}=f(x^2)$, and that
$g(x)=\sum_{n\ge 0}{a_{2n+1}x^{n}}=f(x)^2$ implies that $\sum_{n\ge
0}{a_{2n+1}x^{2n+1}}=xf(x^2)^2$. Hence,
\[
\sum_{n\ge 0}{a_nx^n}=f(x^2)+xf(x^2)^2=\frac{1}{1-xf(x^2)}.
\]

\begin{proof}
Let $\pi\in\D^2_{2n}(2143)$. Note that all fixed points must be on
the upper board of $\pi$. Therefore, the lower board of $\pi$ may be
any $213$-avoiding lower board. This accounts for the factor $a_n$.
Now consider the product of two rational functions on the right.
This product corresponds to the fact that the upper board $B$ of
$\pi$ is a concatenation of two objects: the upper board $B'$ of
rows below the lowest (smallest) fixed point, and the upper board
$B''$ of rows not below the lowest fixed point. It is easy to see
that $B'$ may be any $132$-avoiding upper board. Note that $B''$
must necessarily have an even number of rows and that $B''$ is a
concatenation of a sequence of ``slices'' between consecutive fixed
points, where the $i$th slice consists of an even number of rows
below the $(i+1)$-th smallest fixed point but not below the $i$th
smallest fixed point.

Thus, we obtain a block decomposition of the upper board $B$
(similar to the one in the Figure \ref{fig:2143-lower} for lower
boards) into an possibly empty upper board $B'$ and a sequence $B''$
of nonempty upper boards $B''_1,B''_2,\dots$, where each $B''_i$
contains an even number of rows and exactly 1 fixed point of $\pi$.
Taking generating functions yields the product of functions on the
right-hand side of (\ref{eq:2143-fix}).
\end{proof}

\section{Block Decomposition and Dumont permutations avoiding a pair of 4-letter patterns}
The core of the block decomposition approach initiated by Mansour
and Vainshtein lies in the study of the structure of 132-avoiding
permutations, and permutations containing a given number of
occurrences of 132 (see \cite{MV2} and references therein). In this
section, using the block decomposition approach, we consider those
Dumont permutations in $\S_n$ that avoid a pair of patterns of
length $4$ and an arbitrary pattern.

\subsection{$\{1342,1423\}$-avoiding Dumont permutations of the first kind}

Let $A_\tau(x)$ be the generating function for the number of Dumont
permutations of the first kind in $\D_{2n}^1(1342,1423,\tau)$, that
is,
\[
A_\tau(x)=\sum_{n\ge0}|\D_{2n}^1(1342,1423,\tau)|x^n.
\]

We say a permutation $\tau$ is \emph{decreasing-decomposable} (resp.
\emph{increasing-decomposable}) if there exist nonempty
subpermutations $\tau'$ and $\tau''$ such that $\tau=\tau'\tau''$
and each entry of $\tau'$ is bigger (resp. smaller) than each entry
of $\tau''$.

\begin{theorem}\label{thaa}
Let $\tau\in\S_\ell$ be any pattern which is not
decreasing-decomposable with $\tau_i\ne \ell$ for $i=1,\ell-1,\ell$.
Then
\[
A_\tau(x)=s(x)=\frac{1+x-\sqrt{1-6x+x^2}}{2x}.
\]
\end{theorem}
\begin{proof}
By \cite[Theorem 3.4]{Bur}, we have exactly two possibilities for
the block decomposition of an arbitrary Dumont permutation of the
first kind in $\D_{2n}^1(1342,1423)$. Let us write an equation for
$A_\tau(x)$. The contribution of the first decomposition above is
$xA_\tau(x)(A_\tau(x)-1)$. The contribution of the second possible
decomposition is $x(A_\tau(x))^2$. Therefore, by using the three
contributions above we have that
$A_\tau(x)=1+xA_\tau(x)(A_\tau(x)-1)+x(A_\tau(x))^2$, where $1$ is
the contribution of the empty permutation. Solving this equation
gives the desired result.
\end{proof}

Similarly, we have the following results.

\begin{theorem}\label{thab} ~\\[-\baselineskip]
\begin{enumerate}
\item If $\tau'\in\S_{\ell-1}$ is is not
decreasing-decomposable, $\tau'_{\ell-1}\ne \ell-1$, and
$\tau=\tau'\ell$, then
\[
A_\tau(x)=1+\frac{x(A_{\tau'}(x))^2}{1-xA_{\tau'}(x)}.
\]
\item If $\tau=\tau'(\ell-1)\ell\in \S_\ell$ with no restrictions on
$\tau'\S_{\ell-2}$, then
\[
A_\tau(x)=1+\frac{x(A_{\tau'(\ell-1)}(x))^2}{1-xA_{\tau'}(x)}.
\]
\item If $\tau=\tau'\ell(\ell-1)\in \S_\ell$,  with no restrictions on
$\tau'\S_{\ell-2}$, then
\[
A_\tau(x)=\frac{1+x(1-A_{\tau'}(x))-\sqrt{(1+x(1-A_{\tau'}(x)))^2-4x}}{2x}.
\]
\end{enumerate}
\end{theorem}

For example, if $\tau=13245$, then by Theorem~\ref{thab} we have
$A_{13245}(x)=1+\frac{x(A_{1324}(x))^2}{1-xA_{132}(x)}$. Now, using
Theorem~\ref{thab} for $\tau=1324$ we get that
$A_{1324}(x)=1+\frac{x(A_{132}(x))^2}{1-xA_{132}(x)}$, so
\[
A_{13245}(x)=1+\frac{x(1+xA_{132}(x)(A_{132}(x)-1))^2}{(1-xA_{132}(x))^3}.
\]
Finally, using Theorem~\ref{thab} together with $A_1(x)=1$, we get
that $A_{132}(x)=\frac{1-\sqrt{1-4x}}{2x}=C(x)$. Hence, we can use
$C(x)=\frac{1}{1-xC(x)}$ to obtain $A_{13245}(x)=1+(1-x)^2C^3(x)$.
Another interesting example obtained by Theorem~\ref{thab} is
$A_{2143}(x)=C(x)$ (since $A_{21}(x)=1$). In other words,
$|\D_{2n}^1(1342,1423,2143)|=C_n$. In fact, it is easy to see using
block decomposition that $\D_{2n}^1(1342,1423,2143)=\D_{2n}^1(132)$.

\begin{theorem}\label{thae}
Let $\tau'\in\S_{\ell-1}$ be any non-decreasing-decomposable pattern
with $\tau=\ell\tau'$ and $\tau_\ell\ne \ell-1$. Then
\[
A_\tau(x)=\frac{1}{1+x-2xA_\tau'(x)}.
\]
\end{theorem}
\begin{proof}
By \cite[Theorem 3.4]{Bur}, we have exactly two possibilities for
the block decomposition of an arbitrary Dumont permutation of the
first kind in $\D_{2n}^1(1342,1423)$. Let us write an equation for
$A_\tau(x)$. The contribution of the first decomposition above is
$xA_\tau(x)(A_{\tau'}(x)-1)$. The contribution of the second
possible decomposition is $xA_\tau(x)A_{\tau'}(x)$. Therefore, by
using the three contributions above we have that
$A_\tau(x)=1+xA_\tau(x)(A_{\tau'}(x)-1)+xA_\tau(x)A_{\tau'}(x)$,
where $1$ stands for the empty permutation. Solving this equation
gives the desired result.
\end{proof}

Using the above theorems together with $A_1(x)=A_{21}(x)=1$ and
$A_{12}(x)=1+x$ we get
\[
\begin{array}{l|l|l}
\tau & A_\tau(x) & \text{Reference} \\ \hline
 & & \\[-11pt]
1234 & 1+\frac{x^5(x+2)^2}{(1-x)^2(1-x-x^2)} & \text{Theorem~\ref{thab}}\\
1243 & \frac{1}{1-x^2}C\left(\frac{x}{(1-x^2)^2}\right) & \text{Theorem~\ref{thab}}\\
1324 & 1+xC^3(x) & \text{Theorem~\ref{thab}}\\
1342 & s(x) & \text{Theorem~\ref{thaa}}\\
1423 & s(x) & \text{Theorem~\ref{thaa}}\\
1432 & s(x) & \text{Theorem~\ref{thaa}}\\
2134 & 1+\frac{x}{(1-x)^3} & \text{Theorems~\ref{thaa}
and~\ref{thab}.}
\end{array}.
\]

\subsection{$\{2341,2413\}$-avoiding Dumont permutations of the
first kind}

It was noticed in \cite[Theorem~3.5]{Bur} that $\pi\in
\D_{2n}^1(2341,2413)$ if and only if
\begin{itemize}
\item $\pi=(\pi',2n-1,2n,\pi''+2k)$ for $0\le k\le n-2$,
$\pi'\in \D^1_{2k}(2341,2413)$, $\pi''\in
\D^1_{2n-2k-2}(2341,2413)$;

\item $\pi=(\pi',2n,\pi''+2k,2n-1)$ for $0\le k\le n-1$,
$\pi'\in \D^1_{2k}(2341,2413)$, $\pi''\in
\D^1_{2n-2k-2}(2341,2413)$.
\end{itemize}
This representation is called the \emph{block decomposition} of
$\pi\in \D_{2n}^1(2341,2413)$. Let $B_\tau(x)$ be the generating
function for the number of Dumont permutations of the first kind in
$\D_{2n}^1(2341,2413,\tau)$, that is,
$B_\tau(x)=\sum_{n\ge0}|\D_{2n}^1(2341,2413,\tau)|x^n$.

\begin{theorem}\label{thba}
Let $\tau=\ell\tau'\in\S_\ell$ be a pattern with
$\tau_\ell\ne\ell-1$. Then $B_\tau(x)=\frac{1}{1+x-2xB_{\tau'}(x)}$.
\end{theorem}
\begin{proof}
By Theorem \cite[Theorem~3.5]{Bur}, we have exactly two
possibilities for the block decomposition of an arbitrary $\pi\in
\D_{2n}^1(2341,2413)$. Let us write an equation for $B_\tau(x)$. The
contribution of the first decomposition above is
$xB_\tau(x)(B_{\tau'}(x)-1)$. The contribution of the second
possible decomposition is $xB_{\tau}(x)B_{\tau'}(x)$. Therefore,
adding the two cases with the empty permutation we get
\[
B_\tau(x)=xB_\tau(x)(B_{\tau'}(x)-1)+xB_{\tau}(x)B_{\tau'}(x).
\]
Solving this equation we get the desired result.
\end{proof}

Similarly, we have the following result.

\begin{theorem}\label{thbb}
Let $\tau=\ell\tau'(\ell-1)\in\S_\ell$ be a pattern. Then
$$B_\tau(x)=\frac{1+x(1-B_{\tau'}(x))-\sqrt{(1+x(1-B_{\tau'}(x)))^2-4x}}{2x}.$$
\end{theorem}
For example, if $\tau=4123$ or $\tau=312$, then by
Theorem~\ref{thbb} together with $B_1(x)=B_{21}(x)=1$ we have that
$B_\tau(x)=C(x)$.

\emph{Chebyshev polynomials of the second kind}\/  are defined by
$U_r(\cos\theta)=\frac{\sin(r+1)\theta}{\sin\theta}$ for $r\geq0$.
Clearly, $U_r(t)$ satisfies the following recurrence:
\begin{equation}
U_0(t)=1,\ U_1(t)=2t,\ \mbox{and}\ U_r(t)=2tU_{r-1}(t)-U_{r-2}(t)\
\mbox{for all}\ r\geq2.\label{reccheb}
\end{equation}
and, thus, is a polynomial of degree $r$ in $t$ with integer
coefficients. The same recurrence is used to define $U_r(t)$ for
$r<0$ (for example, $U_{-1}(t)=0$ and $U_{-2}(t)=-1$). The following
lemma can be proved by induction and (\ref{reccheb}).

\newcommand\tuv{\left(\frac{u}{2\sqrt{v}}\right)}
\begin{lemma}\label{lbb}
Define $a_m=\frac{1}{u-va_{m-1}}$ for all $m\ge1$, with $a_0=r$.
Then
$$a_m=\frac{U_{m-1}\tuv-r\sqrt{v}U_{m-2}\tuv}{\sqrt{v}\left[U_m\tuv-r\sqrt{v}U_{m-1}\tuv\right]},$$
where $U_m(t)$ is the $m$-th Chebyshev polynomial of the second
kind.
\end{lemma}

\newcommand\ttx{\left(\frac{1+x}{2\sqrt{2x}}\right)}
\begin{corollary}\label{cba}
For any $k\ge0$,
\[
B_{(k+2)\ldots21}(x)=\frac{U_{k-1}\ttx-\sqrt{2x}U_{k-2}\ttx}{\sqrt{2x}\left[U_k\ttx-\sqrt{2x}U_{k-1}\ttx\right]}.
\]
\end{corollary}
\begin{proof}
It is clear that $B_{21}(x)=1$. Hence, Theorem~\ref{thba} together
with Lemma~\ref{lbb} yields the desired result.
\end{proof}

\subsection{$\{1342,2413\}$-avoiding Dumont permutations of the
first kind}

It was noticed in \cite[Theorem~3.6]{Bur} that $\pi\in
\D_{2n}^1(2341,2413)$ if and only if
\begin{itemize}
\item $\pi=(\pi'+2k,2n-1,2n,\pi'')$ for $1\le k\le n-1$, $\pi'\in
\D^1_{2n-2k-2}(1342,2413)$, $\pi''\in \D^1_{2k}(1342,2413)$;

\item $\pi=(\pi',2n,\pi''+2k,2n-1)$ for $0\le k\le n-1$, $\pi'\in
\D^1_{2k}(1342,2413)$, $\pi''\in \D^1_{2n-2k-2}(1342,2413)$.
\end{itemize}
This representation is called the block decomposition of $\pi\in
\D_{2n}^1(1342,2413) $. Let $C_\tau(x)$ be the generating function
for the number of Dumont permutations of the first kind in
$\D_{2n}^1(1342,2413,\tau)$, that is,
$C_\tau(x)=\sum_{n\ge0}|\D_{2n}^1(1342,2413,\tau)|x^n$.

\begin{theorem}\label{thca}
For all $k\ge3$,
\[
C_{12\ldots k}(x)=1+\frac{x}{1-xC_{12\ldots(k-2)}(x)}
\sum_{j=1}^{k-1}(C_{12\ldots j}(x)-C_{12\ldots(j-1)}(x))
C_{12\ldots(k-j)}(x),
\]
with $C_1(x)=1$ and $C_{12}(x)=1+x$.
\end{theorem}
\begin{proof}
It is easy to check the theorem for $k=1,2$, so we can assume
$k\ge3$. As mentioned before, we have exactly two possibilities for
the block decomposition of an arbitrary $\pi\in
\D_{2n}^1(1342,2413)$. Let us write an equation for $C_{12\ldots
k}(x)$. The contribution of the first decomposition above is
$xC_{12\ldots (k-2)}(x)(C_{12\ldots k}(x)-1)$. The contribution of
the second possible decomposition is $x\sum_{j=1}^{k-1}(C_{12\ldots
j}(x)-C_{12\ldots(j-1)}(x))C_{12\ldots(k-j)}(x)$, since if $\pi'$
contains $12\ldots(j-1)$ and avoids $12\ldots j$, then $\pi''$
avoids $12\ldots(k-j)$  (where $j=1,\ldots,k-1$). Therefore, adding
the two cases with the empty permutation we get
\[
C_{12\ldots k}(x)=1+xC_{12\ldots(k-2)}(x)(C_{12\ldots k}(x)-1)
+x\sum_{j=1}^{k-1}(C_{12\ldots
j}(x)-C_{12\ldots(j-1)}(x))C_{12\ldots(k-j)}(x).
\]
Solving this linear equation we get the desired result.
\end{proof}

For example, Theorem~\ref{thca} for $k=3,4$ gives
$C_{123}(x)=\frac{1+2x^2}{1-x}$ and
$C_{1234}(x)=\frac{1-x+x^2+4x^3+x^4}{(1-x)(1-x-x^2)}$.

\begin{theorem}\label{thcb}
For all $k\ge3$,
\[
C_{k\ldots21}(x)=\frac{1+x\sum_{j=2}^{k-1}(C_{j\ldots21}(x)-C_{(j-1)\ldots21}(x))C_{(k+1-j)\ldots21}(x)}
{1-xC_{(k-1)\ldots21}(x)},
\]
with $C_1(x)=C_{21}(x)=1$.
\end{theorem}
\begin{proof}
It is easy to check the theorem for $k=1,2$. Let $k\ge3$. As
mentioned before, we have exactly two possibilities for the block
decomposition of an arbitrary $\pi\in \D_{2n}^1(1342,2413)$. Let us
write an equation for $C_{k\ldots21}(x)$. The contribution of the
first decomposition above is
\[
x\sum_{j=2}^{k-1}C_{(k+1-j)\ldots21}(x)(C_{j\ldots21}(x)-C_{(j-1)\ldots21}(x)),
\]
where $\pi''$ contains $(j-1)\ldots21$ and avoids $j\ldots21$ for
$j=2,\ldots,k-1$. The contribution of the second possible
decomposition is
\[
xC_{k\ldots21}(x)C_{(k-1)\ldots21}(x).
\]
Therefore, adding the two cases with the empty permutation we get
\[
C_{k\ldots21}(x)=1+xC_{(k-1)\ldots21}(x)C_{k\ldots21}(x)+x\sum_{j=2}^{k-1}(C_{j\ldots21}(x)-C_{(j-1)\ldots21}(x))C_{(k+1-j)\ldots21}(x).
\]
Solving this equation we get the desired expression.
\end{proof}

For example, Theorem~\ref{thcb} for $k=3,4$ gives
$C_{321}(x)=\frac{1}{1-x}$ and $C_{4321}(x)=\frac{1-x+x^2}{1-2x}$.


\end{document}